\documentclass[11pt]{article}%
\usepackage{amsfonts}
\usepackage{amsmath}
\usepackage{amssymb}
\usepackage{graphicx}%
\newtheorem{theorem}{Theorem}

\newtheorem{condition}[theorem]{Condition}

\newtheorem{definition}[theorem]{Definition}
\newtheorem{example}[theorem]{Example}

\newtheorem{lemma}[theorem]{Lemma}

\newtheorem{proposition}[theorem]{Proposition}
\newtheorem{remark}[theorem]{Remark}

\newenvironment{proof}[1][Proof]{\noindent\textbf{#1.} }{\ \rule{0.5em}{0.5em}}
\begin{document}

\title{On the reduction of the multidimensional Schr\"{o}dinger equation to \\a first order equation and its relation to the pseudoanalytic \\function theory}
\author{Vladislav V. Kravchenko\\Depto. de Telecomunicaciones, SEPI\\Escuela Superior de Ingenier\'{\i}a Mec\'{a}nica y El\'{e}ctrica\\Instituto Polit\'{e}cnico Nacional\\C.P.07738 M\'{e}xico D.F., \\MEXICO}
\maketitle

\begin{abstract}
Given a particular solution of a one-dimensional stationary Schr\"{o}dinger
equation this equation of second order can be reduced to a first order linear
ordinary differential equation. This is done with the aid of an auxiliary
Riccati differential equation. In the present work \ we show that the same
fact is true in a multidimensional situation also. For simplicity we consider
the case of two or three independent variables. One particular solution of the
Schr\"{o}dinger equation allows us to reduce this second order equation to a
linear first order quaternionic differential equation. As in one-dimensional
case this is done with the aid of an auxiliary quaternionic Riccati equation.
The resulting first order quaternionic equation is equivalent to the static
Maxwell system and is closely related to the Dirac equation. In the case of
two independent variables it is the well known Vekua equation from theory of
pseudoanalytic (or generalized analytic) functions. Nevertheless we show that
even in this case it is very useful to consider not complex valued functions
only, solutions of the Vekua equation but complete quaternionic functions. In
this way the first order quaternionic equation represents two separate Vekua
equations, one of which gives us solutions of the Schr\"{o}dinger equation and
the other one can be considered as an auxiliary equation of a simpler
structure. Moreover for the auxiliary equation we always have the
corresponding Bers generating pair $(F,G)$, the base of the Bers theory of
pseudoanalytic functions, and what is very important, the Bers derivatives of
solutions of the auxiliary equation give us solutions of the main Vekua
equation and as a consequence of the Schr\"{o}dinger equation.  Based on this
fact we obtain an analogue of the Cauchy integral theorem for solutions of the
Schr\"{o}dinger equation. Other results from theory of pseudoanalytic
functions can be written down for solutions of the Schr\"{o}dinger equation.
Moreover, for an ample class of potentials in the Schr\"{o}dinger equation
(which includes for instance all radial potentials), this new approach gives
us a simple procedure allowing to obtain an infinite sequence of solutions of
the Schr\"{o}dinger equation from one known particular solution.

\end{abstract}

\section{Introduction}

Consider the one-dimensional static Schr\"{o}dinger equation
\begin{equation}
u^{\prime\prime}+vu=0\label{schrod1}%
\end{equation}
and the associated Riccati equation
\begin{equation}
y^{\prime}+y^{2}=-v.\label{ricc1}%
\end{equation}
Equation (\ref{schrod1}) is related to the (\ref{ricc1}) by the easily
inverted substitution
\[
y=\frac{u^{\prime}}{u}.
\]
Thus solutions of the Riccati equation (\ref{ricc1}) are simply logarithmic
derivatives of solutions of the Schr\"{o}dinger equation (\ref{schrod1}) and
vice versa solutions of (\ref{schrod1}) are logarithmic antiderivatives of
solutions of (\ref{ricc1}). The generalization of this fact for a
multidimensional situation was obtained in \cite{KKW} (see also \cite{AQA}).
Among the peculiar properties of the Riccati equation stands out an important
theorem of Euler, dating from 1760. \ If a particular solution $y_{0}$ of the
Riccati equation is known, the substitution $y=y_{0}+z$ reduces (\ref{ricc1})
to a Bernoulli equation which in turn is reduced by the substitution
$z=\frac{1}{u}$ to a first order linear equation. \ Thus given a particular
solution of the Riccati equation, it can be linearized and the general
solution can be found in two integrations. As a consequence of this, given a
particular solution of the Schr\"{o}dinger equation (\ref{schrod1}) the
general solution can be found from a first order linear equation. This can be
seen immediately from the factorization of the one-dimensional Schr\"{o}dinger
operator
\begin{equation}
\partial^{2}+v(x)=(\partial+y_{0}(x))(\partial-y_{0}(x))\label{fact1}%
\end{equation}
which is valid if and only if $y_{0}$ is a solution of (\ref{ricc1}).

In the present work we show that given a particular solution of a
multidimensional stationary Schr\"{o}dinger equation this equation of second
order can be reduced to a first order linear quaternionic differential
equation. For doing this we use a quaternionic factorization of the
Schr\"{o}dinger operator proposed in \cite{Swansolo}, \cite{Swan} (see also
\cite{KK}) and the results on the quaternionic Riccati equation from
\cite{KKW} and \cite{AQA} where it was shown that having a particular solution
of the quaternionic Riccati equation one can reduce it to a second order
linear equation. Here we show that the similarity with the one-dimensional
situation is much closer, and one particular solution is sufficient to reduce
the quaternionic Riccati equation to a first order linear equation. The
resulting first order quaternionic equation is equivalent to the static
Maxwell system and is closely related to the Dirac equation. In the case of
two independent variables it is the well known Vekua equation from theory of
pseudoanalytic (or generalized analytic) functions (see, e.g., \cite{Begehr},
\cite{Berskniga}, \cite{BersStat}, \cite{Tutschke}, \cite{Vekua}). We show
that even in this case it is very useful to consider not only complex valued
functions, solutions of the Vekua equation but complete quaternionic
functions. In this way the first order quaternionic equation represents two
separate Vekua equations, one of which gives us solutions of the
Schr\"{o}dinger equation and the other one can be considered as an auxiliary
equation of a simpler structure. Moreover for the auxiliary equation we always
have in explicit form the corresponding Bers generating pair $(F,G)$, the base
of Bers' theory of pseudoanalytic functions, and what is very important, the
Bers derivatives of solutions of the auxiliary equation give us solutions of
the main Vekua equation and as a consequence of the Schr\"{o}dinger equation.
Based on this fact, for example, we obtain an analogue of the Cauchy integral
theorem for solutions of the Schr\"{o}dinger equation. Other results from
theory of pseudoanalytic functions can be written down for solutions of the
Schr\"{o}dinger equation. Moreover, for an ample class of potentials in the
Schr\"{o}dinger equation (which includes for instance all radial potentials),
this new approach gives us a simple procedure allowing to obtain an infinite
sequence of solutions of the Schr\"{o}dinger equation from one known
particular solution.

Besides this introduction the paper contains four sections. In Section 2 we
introduce necessary notations from quaternionic analysis. In Section 3 we
prove a spatial generalization of the Euler theorem for the Riccati equation
and show how a particular solution of the Schr\"{o}dinger equation allows us
to reduce it to a first order quaternionic equation. We observe that in the
case of two independent variables this first order equation represents two
separate Vekua equations. In order to apply theory of pseudoanalytic functions
to the resulting Vekua equations, in Section 4 we introduce some necessary
definitions and results from Bers' theory. Finally in Section 5 we show how
all the machinery of this quite forgotten mathematical theory allows us to
obtain surprising results for the Schr\"{o}dinger equation starting with an
analogue of the Cauchy integral theorem and including infinite sequences of
solutions generated by one particular solution.

\section{Notations from quaternionic analysis}

We will consider the algebra $\mathbb{H}(\mathbb{C})$ of complex quaternions
or biquaternions which have the form $q=q_{0}+$ $q_{1}\mathbf{i}%
+q_{2}\mathbf{j}+q_{3}\mathbf{k},$ where $\{q_{k}\}\subset\mathbb{C}$, and
$\mathbf{i}$, $\mathbf{j}$, $\mathbf{k}$ are the quaternionic imaginary units. 

The vectorial representation of a complex quaternion will be used. Namely,
each complex quaternion $q$ is a sum of a scalar $q_{0}$ and of a vector
$\mathbf{q}$:
\[
q=\operatorname*{Sc}(q)+\operatorname*{Vec}(q)=q_{0}+\mathbf{q},
\]
where $\mathbf{q}=q_{1}\mathbf{i}+q_{2}\mathbf{j}+q_{3}\mathbf{k}$. The purely
vectorial complex quaternions ($\operatorname*{Sc}(q)=0$) are identified with
vectors from $\mathbb{C}^{3}$. Note that $\mathbf{q}^{2}=-<\mathbf{q}%
,\mathbf{q}>$ where $<\mathbf{\cdot},\mathbf{\cdot}>$ denotes the usual scalar product.

By $M^{p}$ we denote the operator of multiplication by a complex quaternion
$p$ from the right-hand side: $M^{p}q=q\cdot p$.

More information on the structure of the algebra of complex quaternions can be
found for example in \cite{AQA} or \cite{KSbook}.

Let $q$ be a complex quaternion valued differentiable function of
$\mathbf{x}=(x_{1},x_{2},x_{3})$. Denote
\[
Dq=\mathbf{i}\frac{\partial}{\partial x_{1}}q+\mathbf{j}\frac{\partial
}{\partial x_{2}}q+\mathbf{k}\frac{\partial}{\partial x_{3}}q.
\]
This expression can be rewritten in vector form as follows%

\[
Dq=-\operatorname*{div}\mathbf{q}+\operatorname*{grad}q_{0}%
+\operatorname*{rot}\mathbf{q}.
\]
That is, $\operatorname*{Sc}(Dq)=-\operatorname*{div}\mathbf{q}$ and
$\operatorname*{Vec}(Dq)=\operatorname*{grad}q_{0}+\operatorname*{rot}%
\mathbf{q}$. Let us notice that $D^{2}=-\Delta$. 

If $q_{0}$ is a scalar function then $Dq_{0}$ coincides with
$\operatorname*{grad}q_{0}$. The expression $Dq_{0}/q_{0}$ will be called the
logarithmic derivative of $q_{0}$.

\section{Reduction of the Schr\"{o}dinger equation to a first order
quaternionic equation}

Consider the equation
\begin{equation}
(-\Delta+u)f=0\label{Schr}%
\end{equation}
where $\Delta=\frac{\partial^{2}}{\partial x_{1}^{2}}+\frac{\partial^{2}%
}{\partial x_{2}^{2}}+\frac{\partial^{2}}{\partial x_{3}^{2}}$, $f$ and $u$
are complex valued functions. We assume that $f$ is twice differentiable.
Together with (\ref{Schr}) we introduce the following quaternionic equation
\begin{equation}
D\mathbf{q}+\mathbf{q}^{2}=-u\label{Riccati}%
\end{equation}
where $\mathbf{q}$ is a purely vectorial differentiable biquaternion valued function.

\begin{theorem}
\label{ThFact}\cite{Swansolo} For an arbitrary scalar twice differentiable
function $f$ the following equality holds%
\begin{equation}
(D+M^{\mathbf{h}})(D-M^{\mathbf{h}})f=(-\Delta+u)f \label{factQ}%
\end{equation}
if and only if $\mathbf{h}$ is a solution of (\ref{Riccati}).
\end{theorem}

Thus, given a particular solution of (\ref{Riccati}) the Schr\"{o}dinger
operator in (\ref{Schr}) can be factorized.

\begin{theorem}
\label{ThLogder}\cite{KKW} Solutions of (\ref{Schr}) are related to solutions
of (\ref{Riccati}) in the following way. For any nonvanishing solution $f$ of
(\ref{Schr}) its logarithmic derivative
\begin{equation}
\mathbf{q}=\frac{Df}{f} \label{logder}%
\end{equation}
is a solution of (\ref{Riccati}) and any solution $\mathbf{q}$ of
(\ref{Riccati}) is a logarithmic derivative of the form (\ref{logder}) of a
solution of (\ref{Schr}).
\end{theorem}

\begin{proof}
A direct substitution into the equation (\ref{Riccati}) shows us that for a
nonvanishing solution $f$ of (\ref{Schr}) its logarithmic derivative
(\ref{logder}) is a solution of (\ref{Riccati}). Now let us suppose that
$\mathbf{q}$ is a solution of (\ref{Riccati}). From the vector part of
(\ref{Riccati}) we have that $\mathbf{q}$ is a gradient of some scalar
function $\xi$: $\mathbf{q}=\operatorname{grad}\xi$. Then $\mathbf{q}$ can be
represented in the form (\ref{logder}) where $f=e^{\xi}$. Substituting
(\ref{logder}) in (\ref{Riccati}) we obtain that $f$ is a solution of
(\ref{Schr}).
\end{proof}

\begin{remark}
Theorems \ref{ThFact} and \ref{ThLogder} show us that equation (\ref{Riccati})
is a generalization of the Riccati equation. We will call it quaternionic
Riccati equation.
\end{remark}

\begin{lemma}
\label{LemmaMax}\cite{KrZAA02} For a nonvanishing scalar differentiable
function $\varepsilon$ there exists a one-to-one correspondence between
solutions of the static Maxwell system
\begin{equation}
\operatorname{div}(\varepsilon\mathbf{E})=0, \label{M1}%
\end{equation}%
\begin{equation}
\operatorname{rot}\mathbf{E}=0, \label{M2}%
\end{equation}
and solutions of the equation
\begin{equation}
(D+M^{\mathbf{h}})\mathbf{F}=0 \label{D+Mh}%
\end{equation}
where
\[
\mathbf{h}=\frac{D\sqrt{\varepsilon}}{\sqrt{\varepsilon}}.
\]
Vector $\mathbf{E}$ is a solution of (\ref{M1}), (\ref{M2}) if and only if the
vector $\mathbf{F}=\sqrt{\varepsilon}\mathbf{E}$ is a solution of (\ref{D+Mh}).
\end{lemma}

\begin{proof}
The system (\ref{M1}), (\ref{M2}) can be rewritten in the form%
\[
D\mathbf{E}=<\frac{\operatorname{grad}\varepsilon}{\varepsilon},\mathbf{E}>.
\]
Let us make a simple observation: the scalar product of two vectors
$\mathbf{p}$ and $\mathbf{q}$ can be written as follows%
\[
<\mathbf{p},\mathbf{q}>=-\frac{1}{2}(^{\mathbf{p}}M+M^{\mathbf{p}})\mathbf{q}.
\]
Then we have
\begin{equation}
(D+\frac{1}{2}\frac{\operatorname{grad}\varepsilon}{\varepsilon}%
)\mathbf{E}=-\frac{1}{2}M^{\frac{\operatorname{grad}\varepsilon}{\varepsilon}%
}\mathbf{E.}\label{Vsp1}%
\end{equation}
Note that
\[
\frac{1}{2}\frac{\operatorname{grad}\varepsilon}{\varepsilon}=\frac
{\operatorname{grad}\sqrt{\varepsilon}}{\sqrt{\varepsilon}}.
\]
Then equation (\ref{Vsp1}) can be rewritten in the following form%
\begin{equation}
\frac{1}{\sqrt{\varepsilon}}D(\sqrt{\varepsilon}\mathbf{E)}+\mathbf{Eh}%
=0\label{Min31}%
\end{equation}
where $\mathbf{h}=D\sqrt{\varepsilon}/\sqrt{\varepsilon}$. Introducing the
notation $\mathbf{F}=\sqrt{\varepsilon}\mathbf{E}$ and multiplying
(\ref{Min31}) by $\sqrt{\varepsilon}$ we obtain the equivalence of the system
(\ref{M1})-(\ref{M2}) to (\ref{D+Mh}).
\end{proof}

\begin{lemma}
\label{EulerOld}\cite{KKW} Let $\mathbf{h}$ be an arbitrary particular
solution of (\ref{Riccati}) (then as was mentioned above it is a gradient of
some scalar function $\xi$). \ The general solution of (\ref{Riccati}) has the
form
\begin{equation}
\mathbf{q}=\mathbf{h}+\mathbf{g},\label{sum}%
\end{equation}
where $\mathbf{g}=(\operatorname{grad}\Psi)/\Psi$ and $\Psi$ is a general
solution of the equation
\begin{equation}
\triangle\Psi+2\left\langle \operatorname{grad}\xi,\operatorname{grad}%
\Psi\right\rangle =0,\label{trans}%
\end{equation}
or equivalently of
\begin{equation}
\operatorname*{div}(e^{2\xi}\operatorname{grad}\Psi)=0.\label{div}%
\end{equation}

\end{lemma}

\begin{proof}
Substituting (\ref{sum}) in (\ref{Riccati}) gives
\begin{equation}
D\mathbf{g}-2\left\langle \mathbf{h},\mathbf{g}\right\rangle +\mathbf{g}%
^{2}=0.\label{ours}%
\end{equation}
Note that the vector part of (\ref{ours}) is $\operatorname*{rot}%
\mathbf{g}=0,$ so that
\[
\mathbf{g}=\operatorname{grad}\Phi
\]
for some function $\Phi.$ \ If $\Psi=e^{\Phi},$ this is equivalent to
\[
\mathbf{g}=(\operatorname{grad}\Psi)/\Psi.
\]
Equation (\ref{ours}), written in terms of $\Psi,$ is
\[
-\frac{1}{\Psi^{2}}(\operatorname{grad}\Psi)^{2}-\frac{1}{\Psi}\triangle
\Psi-\frac{2}{\Psi}\left\langle \operatorname{grad}\xi,\operatorname{grad}%
\Psi\right\rangle +\frac{1}{\Psi^{2}}(\operatorname{grad}\Psi)^{2}=0,
\]
so that (\ref{ours}) is equivalent to
\[
\triangle\Psi+2\left\langle \operatorname{grad}\xi,\operatorname{grad}%
\Psi\right\rangle =0.
\]
Noting that
\[
\operatorname*{div}(e^{2\xi}\operatorname{grad}\Psi)=2e^{2\xi}\left\langle
\operatorname{grad}\xi,\operatorname{grad}\Psi\right\rangle +e^{2\xi}%
\triangle\Psi=e^{2\xi}(\triangle\Psi+2\left\langle \operatorname{grad}%
\xi,\operatorname{grad}\Psi\right\rangle ),
\]
this equation can be rewritten in the form
\[
\operatorname*{div}(e^{2\xi}\operatorname{grad}\Psi)=0.
\]

\end{proof}

Now we are ready to prove a generalization of the Euler theorem for the
quaternionic Riccati equation.

\begin{theorem}
\label{ThEuler}(Euler's theorem for the quaternionic Riccati equation) Let
$\mathbf{h}=\operatorname{grad}\xi$ be a particular solution of (\ref{Riccati}%
). The general solution of the quaternionic Riccati equation has the form
$\mathbf{q}=\mathbf{h}+\mathbf{g}$ where $\mathbf{g}=\frac{D\Psi}{\Psi}$ and
$\Psi$ is obtained from the equation
\begin{equation}
\operatorname{grad}\Psi=e^{-\xi}\mathbf{F} \label{PsiF}%
\end{equation}
where $\mathbf{F}$ is the general solution of (\ref{D+Mh}).
\end{theorem}

\begin{proof}
According to Lemma \ref{EulerOld} it is sufficient to prove that $\Psi$ is a
solution of (\ref{trans}) (or what is the same of (\ref{div})) if and only if
the vector $\mathbf{F}=e^{\xi}\operatorname{grad}\Psi$ is a solution of
(\ref{D+Mh}). Let us notice that if $\Psi$ is a solution of (\ref{trans}) then
the vector $\mathbf{E}=\operatorname{grad}\Psi$ is a solution of the system
(\ref{M1})-(\ref{M2}) where $\varepsilon=e^{2\xi}$ and vice versa if
$\mathbf{E}$ is a solution of (\ref{M1})-(\ref{M2}) then it is a gradient of
some function $\Psi$ which is necessarily a solution of (\ref{trans}). Now due
to Lemma \ref{LemmaMax}, $\operatorname{grad}\Psi$ is a solution of
(\ref{M1})-(\ref{M2}) if and only if $\mathbf{F}=e^{\xi}\operatorname{grad}%
\Psi$ is a solution of the equation (\ref{D+Mh}) where $\mathbf{h}%
=\frac{De^{\xi}}{e^{\xi}}=\operatorname{grad}\xi$.
\end{proof}

Thus, given a particular solution of the quaternionic Riccati equation, the
general solution reduces to the linear first order equation (\ref{D+Mh}),
exactly as in one-dimensional situation.

Using Theorem \ref{ThLogder} we immediately arrive at the following result for
the Schr\"{o}dinger equation.

\begin{theorem}
\label{ThSchr}Let $f_{0}$ be a nonvanishing particular solution of
(\ref{Schr}) and $\mathbf{F}$ be the general solution of (\ref{D+Mh}) where
$\mathbf{h}=Df_{0}/f_{0}$. Then the general solution $f$ of (\ref{Schr}) is
the logarithmic antiderivative of $\mathbf{q}$: $\frac{Df}{f}=\mathbf{q}$,
where $\mathbf{q}=\mathbf{h}+\mathbf{g}$, $\mathbf{g}=\frac{D\Psi}{\Psi}$ and
$\Psi$ is obtained from the equation
\begin{equation}
\operatorname{grad}\Psi=\frac{\mathbf{F}}{f_{0}}. \label{PsiFf}%
\end{equation}

\end{theorem}

\begin{proof}
From Theorem \ref{ThLogder} we have that $f$ is the logarithmic antiderivative
of $\mathbf{q}$ (equation (\ref{logder})), where $\mathbf{q}$ is the general
solution of (\ref{Riccati}). For a particular solution $f_{0}$ of (\ref{Schr})
the vector
\[
\mathbf{h}=\frac{Df_{0}}{f_{0}}=\operatorname{grad}\ln f_{0}%
\]
is a particular solution of (\ref{Riccati}). Due to Theorem \ref{ThEuler} the
general solution of (\ref{Riccati}) has the form $\mathbf{q}=\mathbf{h}%
+\mathbf{g}$, where $\mathbf{g}=\frac{D\Psi}{\Psi}$ and $\Psi$ is obtained
from the equation
\[
\operatorname{grad}\Psi=e^{-\ln f_{0}}\mathbf{F=}\frac{\mathbf{F}}{f_{0}}%
\]
where $\mathbf{F}$ is the general solution of (\ref{D+Mh}).
\end{proof}

\begin{remark}
\label{RemLog}The logarithmic antiderivative of $\mathbf{q}$ always exists and
can be obtained easily. Being a solution of (\ref{Riccati}) $\mathbf{q}$ is
necessarily a gradient of some function $\Phi$ which can be constructed
analytically. Then $f$ has the form $f=Ce^{\Phi}$, where $C$ is a complex constant.
\end{remark}

\begin{remark}
From Theorem \ref{ThSchr} it follows that for any vector $\mathbf{F,}$
solution of (\ref{D+Mh}) with $\mathbf{h}=Df_{0}/f_{0}$, the vector
$\mathbf{F}/f_{0}$ must be a gradient of some scalar function $\Psi$, and this
is really true. Let us show that indeed $\operatorname{rot}(\mathbf{F}%
/f_{0})=0$. Note that this condition is equivalent to the equality
$\operatorname{Vec}(D(\mathbf{F}/f_{0}))=0$. Consider
\[
D(\frac{\mathbf{F}}{f_{0}})=\frac{1}{f_{0}}D\mathbf{F}-\frac{Df_{0}}{f_{0}%
^{2}}\mathbf{F}=-(\mathbf{F}\frac{Df_{0}}{f_{0}^{2}}+\frac{Df_{0}}{f_{0}^{2}%
}\mathbf{F})=2<\mathbf{F},\frac{Df_{0}}{f_{0}^{2}}>.
\]
Thus $\operatorname{Vec}(D(\mathbf{F}/f_{0}))=0$ and hence the vector
$\mathbf{F}/f_{0}$ is a gradient.
\end{remark}

\begin{remark}
\label{RemChain}Let us summarize the results of this section in the following
chain of actions which one should follow in order to obtain solutions of
(\ref{Schr}) from solutions of (\ref{D+Mh}).

Given a nonvanishing particular solution $f_{0}$ of the Schr\"{o}dinger
equation (\ref{Schr}) we construct the vector $\mathbf{h}=Df_{0}/f_{0}$ and
consider equation (\ref{D+Mh}). Taking a solution $\mathbf{F}$ of (\ref{D+Mh})
we find $\Psi$ from (\ref{PsiFf}). Then we construct the vectors
$\mathbf{g}=D\Psi/\Psi$ and $\mathbf{q}=\mathbf{h}+\mathbf{g}$. $\mathbf{q}$
is necessarily a gradient of some scalar function $\Phi$. Finding $\Phi$ we
finally obtain a solution of (\ref{Schr}) as $f=Ce^{\Phi}$, where $C$ is a
complex constant.

Thus given a particular solution of the Schr\"{o}dinger equation (\ref{Schr}),
the general solution reduces to the first order equation (\ref{D+Mh}). It is
interesting to note that due to Lemma \ref{LemmaMax}, equation (\ref{D+Mh}) is
equivalent to the static Maxwell system. It is closely related (see \cite{KK}
and \cite{KrBelt}) also to the Dirac equation as well as to the Beltrami
fields which are solutions of the equation $\operatorname{rot}\mathbf{f}%
+\alpha\mathbf{f}=0$ (see, e.g., \cite{Stratis} and \cite{Lak}).
\end{remark}

The following statement gives us the way to transform solutions of the
Schr\"{o}dinger equation into solutions of (\ref{D+Mh}).

\begin{proposition}
\label{PropSchrVek}Let $f_{1}$ be another nonvanishing solution of the
Schr\"{o}dinger equation (\ref{Schr}). Then the ratio $\Psi=f_{1}/f_{0}$ is a
solution of the equation%
\begin{equation}
\operatorname*{div}(f_{0}^{2}\operatorname*{grad}\Psi)=0 \label{divf0}%
\end{equation}
and the vector $\mathbf{F}=f_{0}D(f_{1}/f_{0})$ is a solution of equation
(\ref{D+Mh}) where $\mathbf{h}=Df_{0}/f_{0}$.
\end{proposition}

\begin{proof}
By Theorem \ref{ThLogder} the vector $\mathbf{q}=Df_{1}/f_{1}$ is a solution
of (\ref{Riccati}). Then consider the vector $\mathbf{g}$ from Lemma
\ref{EulerOld}:
\[
\mathbf{g}=\mathbf{q}-\mathbf{h}=\frac{Df_{1}}{f_{1}}-\frac{Df_{0}}{f_{0}%
}=D(f_{1}/f_{0})/(f_{1}/f_{0}).
\]
Thus we have that $\mathbf{g}=D\Psi/\Psi$ where $\Psi=f_{1}/f_{0},$ and by
Lemma \ref{EulerOld}, $\Psi$ satisfies (\ref{divf0}).

From (\ref{PsiFf}) we obtain that $\mathbf{F}=f_{0}D\Psi=f_{0}D(f_{1}/f_{0})$
is a solution of (\ref{D+Mh}).
\end{proof}

Let us consider equation (\ref{D+Mh}) for a biquaternion valued function $p$
whose scalar part is not necessarily zero%
\begin{equation}
(D+M^{\mathbf{h}})p=0\label{D+Mhp}%
\end{equation}
and let us use the following representation for biquaternions:%
\[
p=P_{1}+P_{2}\mathbf{j},
\]
where $P_{1}=p_{0}+p_{3}\mathbf{k}$ and $P_{2}=p_{2}-p_{1}\mathbf{k}$. Then
$D=D_{1}+D_{2}\mathbf{j}$, where $D_{1}=\partial_{3}\mathbf{k}$,
$D_{2}=\partial_{2}-\partial_{1}\mathbf{k}$ and $\mathbf{h}=H_{1}%
+H_{2}\mathbf{j}$, where $H_{1}=h_{3}\mathbf{k}$, $H_{2}=h_{2}-h_{1}%
\mathbf{k}$. Using these notations, equation (\ref{D+Mhp}) can be rewritten as
the following system%
\begin{equation}
D_{1}P_{1}-D_{2}\overline{P}_{2}+H_{1}P_{1}-\overline{H}_{2}P_{2}%
=0,\label{sys1}%
\end{equation}%
\begin{equation}
D_{2}\overline{P}_{1}+D_{1}P_{2}+H_{2}P_{1}+\overline{H}_{1}P_{2}%
=0.\label{sys2}%
\end{equation}
Now let us suppose that both $p$ and $\mathbf{h}$ do not depend on $x_{3}$.
Then system (\ref{sys1}), (\ref{sys2}) turns into the pair of decoupled
equations:%
\begin{equation}
\overline{D}_{2}P_{1}+\overline{H_{2}P_{1}}=0\label{VekP1}%
\end{equation}
and%
\begin{equation}
\overline{D}_{2}P_{2}+H_{2}\overline{P}_{2}=0\label{VekP2}%
\end{equation}
which are nothing but Vekua's equations describing pseudoanalytic or
generalized analytic functions (see, e.g., \cite{Berskniga} and \cite{Vekua}).

\begin{remark}
Here we should mention that in general the components $p_{0}$,...,$p_{3}$ as
well as $h_{1}$, $h_{2}$ can be complex valued functions, hence $P_{1}$,
$P_{2}$ and $H_{2}$ can be bicomplex. Nevertheless this detail is
insignificant for what follows, because all results from Bers' theory which
will be used in the subsequent sections are valid for bicomplex solutions
also. Of course, when $u$ in (\ref{Schr}) is a real valued function we can
consider real valued solutions of (\ref{Schr}) only. In that case
(\ref{VekP1}) and (\ref{VekP2}) are usual Vekua's equations.
\end{remark}

\begin{remark}
\label{RemFG}In what follows we consider $u$ and $f$ in (\ref{Schr}) being
independent of $x_{3}$. Then given a particular solution $f_{0}$ of
(\ref{Schr}), the general solution reduces to equation (\ref{VekP2}) (which in
this case is equivalent to (\ref{D+Mh})). Thus we are primarily interested in
solutions of (\ref{VekP2}), and (\ref{VekP1}) can be considered as an
auxiliary equation. Nevertheless as we will see in Section 5 equations
(\ref{VekP1}) and (\ref{VekP2}) are closely related to each other. With the
aid of Bers' theory solutions of (\ref{VekP2}) can be obtained from solutions
of (\ref{VekP1}) and it is interesting to notice that by construction we
always have at least two solutions of (\ref{VekP1}) in explicit form. It is
easy to see that the functions
\[
F=\frac{1}{f_{0}}\qquad\text{and}\qquad G=f_{0}\mathbf{k}%
\]
are solutions of (\ref{D+Mhp}) where $\mathbf{h}=Df_{0}/f_{0}$ and
consequently they are solutions of (\ref{VekP1}).
\end{remark}

\section{Some definitions and results from Bers' theory}

Bers' theory of pseudoanalytic functions was essentially developed in
\cite{Berskniga} (see also \cite{BersStat}). It is based on the so called
generating pair, a pair of complex functions $F$ and $G$ satisfying the
inequality
\begin{equation}
\operatorname*{Im}(\overline{F}G)>0 \label{ImFG}%
\end{equation}
in some domain of interest $\Omega$ which may coincide with the whole complex
plane. $F$ and $G$ are assumed to possess partial derivatives with respect to
the real variables $x$ and $y$. In this case the operators $\partial
_{\overline{z}}=\frac{\partial}{\partial x}+i\frac{\partial}{\partial y}$ and
$\partial_{z}=\frac{\partial}{\partial x}-i\frac{\partial}{\partial y}$ can be
applied (usually these operators are introduced with the factor $1/2$,
nevertheless here it is somewhat more convenient to consider them without it)
and the following characteristic coefficients of the pair $(F,G)$ can be
defined%
\[
a_{(F,G)}=-\frac{\overline{F}G_{\overline{z}}-F_{\overline{z}}\overline{G}%
}{F\overline{G}-\overline{F}G},\qquad b_{(F,G)}=\frac{FG_{\overline{z}%
}-F_{\overline{z}}G}{F\overline{G}-\overline{F}G},
\]

\[
A_{(F,G)}=-\frac{\overline{F}G_{z}-F_{z}\overline{G}}{F\overline{G}%
-\overline{F}G},\qquad B_{(F,G)}=\frac{FG_{z}-F_{z}G}{F\overline{G}%
-\overline{F}G},
\]
where the subindex $\overline{z}$ or $z$ means the application of
$\partial_{\overline{z}}$ or $\partial_{z}$ respectively.

Every complex function $w$ defined in a subdomain of $\Omega$ admits the
unique representation $w=\phi F+\psi G$ where the functions $\phi$ and $\psi$
are real valued. Sometimes it is convenient to associate with the function $w$
the function $\omega=\phi+i\psi$. The correspondence between $w$ and $\omega$
is one-to-one.

Bers introduces the notion of the $(F,G)$-derivative of a function $w$ which
exists and has the form
\begin{equation}
\overset{\cdot}{w}=\phi_{z}F+\psi_{z}G=w_{z}-A_{(F,G)}w-B_{(F,G)}\overline{w}
\label{FGder}%
\end{equation}
if and only if
\begin{equation}
\phi_{\overline{z}}F+\psi_{\overline{z}}G=0. \label{phiFpsiG}%
\end{equation}
This last equation can be rewritten in the following form%
\[
w_{\overline{z}}=a_{(F,G)}w+b_{(F,G)}\overline{w}%
\]
which we call the Vekua equation. Solutions of this equation are called
$(F,G)$-pseudoanalytic functions. If $w$ is $(F,G)$-pseudoanalytic, the
associated function $\omega$ is called $(F,G)$-pseudoanalytic of second kind.

\begin{remark}
The functions $F$ and $G$ are $(F,G)$-pseudoanalytic, and $\overset{\cdot}%
{F}\equiv\overset{\cdot}{G}\equiv0$.
\end{remark}

\begin{definition}
\label{DefSuccessor}Let $(F,G)$ and $(F_{1},G_{1})$ - be two generating pairs
in $\Omega$. $(F_{1},G_{1})$ is called \ successor of $(F,G)$ and $(F,G)$ is
called predecessor of $(F_{1},G_{1})$ if%
\[
a_{(F_{1},G_{1})}=a_{(F,G)}\qquad\text{and}\qquad b_{(F_{1},G_{1})}%
=-B_{(F,G)}\text{.}%
\]

\end{definition}

The importance of this definition becomes obvious from the following statement.

\begin{theorem}
\label{ThBersDer}Let $w$ be an $(F,G)$-pseudoanalytic function and let
$(F_{1},G_{1})$ be a successor of $(F,G)$. Then $\overset{\cdot}{w}$ is an
$(F_{1},G_{1})$-pseudoanalytic function.
\end{theorem}

\begin{definition}
\label{DefAdjoint}Let $(F,G)$ be a generating pair. Its adjoint generating
pair $(F,G)^{\ast}=(F^{\ast},G^{\ast})$ is defined by the formulas%
\[
F^{\ast}=-\frac{2\overline{F}}{F\overline{G}-\overline{F}G},\qquad G^{\ast
}=\frac{2\overline{G}}{F\overline{G}-\overline{F}G}.
\]

\end{definition}

\begin{theorem}%
\[
(F,G)^{\ast\ast}=(F,G),
\]%
\[
a_{(F^{\ast},G^{\ast})}=-a_{(F,G)},\qquad A_{(F^{\ast},G^{\ast})}=-A_{(F,G)},
\]%
\[
b_{(F^{\ast},G^{\ast})}=-\overline{B}_{(F,G)},\qquad B_{(F^{\ast},G^{\ast}%
)}=-\overline{b}_{(F,G)}.
\]

\end{theorem}

\begin{lemma}
If $(F_{1},G_{1})$ is a successor of $(F,G)$ then $(F,G)^{\ast}$ is a
successor of $(F_{1},G_{1})^{\ast}$.
\end{lemma}

The $(F,G)$ integral of $w$ on a rectifiable curve $\Gamma$ is, by definition,%
\[
\int_{\Gamma}wd_{(F,G)}z=\operatorname*{Re}\int_{\Gamma}F^{\ast}%
wdz-i\operatorname*{Re}\int_{\Gamma}G^{\ast}wdz.
\]
Another important integral is also needed%
\[
\ast\int_{\Gamma}wd_{(F,G)}z=\operatorname*{Re}\int_{\Gamma}G^{\ast
}wdz+i\operatorname*{Re}\int_{\Gamma}F^{\ast}wdz
\]
(we follow the notations of L. Bers).

A continuous function $w$ defined in a domain $\Omega$ is called
$(F,G)$-integrable if for every closed curve $\Gamma$ situated in a simply
connected subdomain of $\Omega$,
\[
\int_{\Gamma}wd_{(F,G)}z=0.
\]

\begin{theorem}
\label{ThCauchyPa}An $(F,G)$-derivative $\overset{\cdot}{w}$ of an
$(F,G)$-pseudoanalytic function $w$ is $(F,G)$-integrable and $\ast\int
_{z_{0}}^{z_{1}}\overset{\cdot}{w}d_{(F,G)}z=\omega(z_{1})-\omega(z_{0})$.
\end{theorem}

The integral $\ast\int_{z_{0}}^{z_{1}}\overset{\cdot}{w}d_{(F,G)}z$ is called
$(F,G)$-antiderivative of $\overset{\cdot}{w}$.

\begin{theorem}
\label{ThMoreraPa}Let $(F,G)$ be a predecessor of $(F_{1},G_{1})$. A
continuous function is $(F_{1},G_{1})$-pseudoanalytic if and only if it is
$(F,G)$-integrable.
\end{theorem}

\section{Applications of Bers' theory to the Schr\"{o}dinger equation}

Let us return to equations (\ref{VekP1}) and (\ref{VekP2}) which in a two
dimensional case are equivalent to the quaternionic equation (\ref{D+Mhp}). In
order to use Bers' notations from the preceding section we rewrite
(\ref{VekP1}) and (\ref{VekP2}) in the following form%

\begin{equation}
w_{\overline{z}}=b\overline{w} \label{Vek1}%
\end{equation}
and%
\begin{equation}
v_{\overline{z}}=\overline{b}\overline{v} \label{Vek2}%
\end{equation}
where $z=x+iy$, $x=x_{2}$, $y=x_{1}$ and instead of the imaginary unit
$\mathbf{k}$ we write $i$. It is easy to see that $b=-\overline{H}%
_{2}=-\partial_{\overline{z}}f_{0}/f_{0}$ and\ $w=P_{1}$, $v=P_{2}$.

As was mentioned above (Remark \ref{RemFG}) for equation (\ref{Vek1}) we know
always two solutions%

\begin{equation}
F=\frac{1}{f_{0}}\qquad\text{and}\qquad G=if_{0}, \label{Berspair}%
\end{equation}
which obviously fulfill (\ref{ImFG}). Thus $(F,G)$ is a generating pair
corresponding to equation (\ref{Vek1}). We have
\[
a_{(F,G)}=0,\qquad b_{(F,G)}=b=-\frac{\partial_{\overline{z}}f_{0}}{f_{0}},
\]

\[
A_{(F,G)}=0,\qquad B_{(F,G)}=-\frac{\partial_{z}f_{0}}{f_{0}}.
\]
According to Definition \ref{DefSuccessor} the characteristic coefficients for
a successor of $(F,G)$ have the form
\[
a_{(F_{1},G_{1})}=0,\qquad b_{(F_{1},G_{1})}=\frac{\partial_{z}f_{0}}{f_{0}%
}=-\overline{b}.
\]
Then due to Theorem \ref{ThBersDer}, if $w$ is a solution of (\ref{Vek1}) then
its $(F,G)$-derivative is a solution of the equation
\begin{equation}
W_{\overline{z}}=-\overline{b}\overline{W}, \label{WVek}%
\end{equation}
but solutions of the last equation multiplied by $i$ become solutions of
(\ref{Vek2}) and vice versa. Thus we obtain the following statement.

\begin{theorem}
Let $w$ be a solution of (\ref{Vek1}). Then the function
\[
v=i\overset{\cdot}{w}=i(w_{z}+\frac{\partial_{z}f_{0}}{f_{0}}\overline{w})
\]
is a solution of (\ref{Vek2}).
\end{theorem}

It is easy to see that according to Definition \ref{DefAdjoint}:
\[
F^{\ast}=-\frac{i}{f_{0}},\qquad G^{\ast}=f_{0}%
\]
and%
\[
b_{(F^{\ast},G^{\ast})}=-\overline{B}_{(F,G)}=-b.
\]
Thus the $(F,G)$ integral of a function $W$ is defined as follows%
\[
\int_{\Gamma}Wd_{(F,G)}z=-\operatorname*{Re}\int_{\Gamma}\frac{i}{f_{0}%
}Wdz-i\operatorname*{Re}\int_{\Gamma}f_{0}Wdz
\]%
\[
=\operatorname*{Im}\int_{\Gamma}\frac{W}{f_{0}}dz-i\operatorname*{Re}%
\int_{\Gamma}f_{0}Wdz.
\]
From Theorems \ref{ThCauchyPa} and \ref{ThMoreraPa} we obtain the following result.

\begin{theorem}
Let $v$ be a solution of (\ref{Vek2}) in a domain $\Omega$. Then for every
closed curve $\Gamma$ situated in a simply connected subdomain of $\Omega$,%
\begin{equation}
\operatorname*{Re}\int_{\Gamma}\frac{v}{f_{0}}dz+i\operatorname*{Im}%
\int_{\Gamma}f_{0}vdz=0. \label{CauchyTh}%
\end{equation}

\end{theorem}

\begin{proof}
For any solution $v$ of (\ref{Vek2}) the function $W=iv$ is a solution of
(\ref{WVek}). As (\ref{WVek}) corresponds to a successor of $(F,G)$, by
Theorem \ref{ThCauchyPa} $W$ is $(F,G)$-integrable. That is%
\[
\operatorname*{Im}\int_{\Gamma}\frac{W}{f_{0}}dz-i\operatorname*{Re}%
\int_{\Gamma}f_{0}Wdz=0.
\]
Now substituting $iv$ instead of $W$ we obtain (\ref{CauchyTh}).
\end{proof}

In order to analyze the meaning of this result for solutions of the
Schr\"{o}dinger equation let us rewrite some statements from Section 3 in our
\textquotedblleft two-dimensional\textquotedblright\ notations.

Consider the equation%
\begin{equation}
(-\partial_{z}\partial_{\overline{z}}+u)f=0 \label{Schrc}%
\end{equation}
where $u$ and $f$ depend on $x$, $y$, and $z=x+iy$. For simplicity we consider
$u$ and $f$ being real valued functions. The corresponding Riccati equation
(\ref{Riccati}) takes the form%
\[
\partial_{\overline{z}}Q+\left\vert Q\right\vert ^{2}=u
\]
where $\mathbf{q}=Q\mathbf{j}$. Equation (\ref{D+Mh}) turns into (\ref{Vek2}),
$\mathbf{F}$ from (\ref{D+Mh}) and $v$ from (\ref{Vek2}) are related by the
equality $\mathbf{F}=v\mathbf{j}$. Then Theorem \ref{ThSchr} in a
two-dimensional situation can be rewritten as follows.

\begin{theorem}
Let $f_{0}$ be a nonvanishing particular solution of (\ref{Schrc}) and $v$ be
the general solution of (\ref{Vek2}) where $b=-\partial_{\overline{z}}%
f_{0}/f_{0}$. Then the general solution $f$ of (\ref{Schrc}) is obtained from
the equation $\partial_{z}f=Qf$, where $Q=\partial_{z}f_{0}/f_{0}+\partial
_{z}\Psi/\Psi$ and $\Psi$ is obtained from the equation $\partial_{z}%
\Psi=v/f_{0}$.
\end{theorem}

As was explained in Remarks \ref{RemLog}-\ref{RemChain}, given a solution $v$
of (\ref{Vek2}), the corresponding solution $f$ of (\ref{Schrc}) can be
constructed analytically. The procedure requires on two steps to reconstruct
the potential function from its gradient.

From Proposition \ref{PropSchrVek} we obtain the following statement.

\begin{proposition}
Let $f_{1}$ be another nonvanishing solution of (\ref{Schrc}). Then the
function
\begin{equation}
v=f_{0}\partial_{z}(f_{1}/f_{0}) \label{vf}%
\end{equation}
is a solution of (\ref{Vek2}), where $b=-\partial_{\overline{z}}f_{0}/f_{0}$.
\end{proposition}

Having this precise relation between solutions of (\ref{Schrc}) and
(\ref{Vek2}) we are able to prove the following result.

\begin{theorem}
\label{ThCauchySchr}(Cauchy's integral theorem for the Schr\"{o}dinger
equation) Let $f_{0}$ and $f_{1}$ be two arbitrary nonvanishing solutions of
(\ref{Schrc}) in a domain $\Omega$. Then for every closed curve $\Gamma$
situated in a simply connected subdomain of $\Omega$,%
\begin{equation}
\operatorname*{Re}\int_{\Gamma}\partial_{z}(\frac{f_{1}}{f_{0}}%
)dz+i\operatorname*{Im}\int_{\Gamma}f_{0}^{2}\partial_{z}(\frac{f_{1}}{f_{0}%
})dz=0. \label{CauchySchrTh}%
\end{equation}

\end{theorem}

\begin{proof}
Substitution of (\ref{vf}) into (\ref{CauchyTh}) gives us the result.
\end{proof}

In Example \ref{ExCauchyTh} we will give a nontrivial example illustrating
this theorem.

From Theorem \ref{ThMoreraPa} and Theorem \ref{ThCauchySchr} we obtain an
analogue of the Morera theorem for the Schr\"{o}dinger equation (\ref{Schrc}).

\begin{theorem}
Let $f_{0}$ be a nonvanishing particular solution of (\ref{Schrc}). The
function $f_{1}$ is a solution of (\ref{Schrc}) also if (\ref{CauchySchrTh})
is valid for every closed curve $\Gamma$ situated in a simply connected
subdomain of $\Omega.$
\end{theorem}

Consider equation (\ref{phiFpsiG}) with the functions (\ref{Berspair}). It
takes the form%
\begin{equation}
\phi_{\overline{z}}+if_{0}^{2}\psi_{\overline{z}}=0. \label{phipsiif2}%
\end{equation}

\begin{proposition}
Let the function $w=\phi F+\psi G$ be $(F,G)$-pseudoanalytic corresponding to
the functions (\ref{Berspair}), that is
\begin{equation}
w=\frac{\phi}{f_{0}}+if_{0}\psi\label{wforVek1}%
\end{equation}
is a solution of (\ref{Vek1}). Then
\[
<\operatorname*{grad}\phi,\operatorname*{grad}\psi>=0.
\]

\end{proposition}

\begin{proof}
The function (\ref{wforVek1}) is a solution of (\ref{Vek1}) if and only if
(\ref{phipsiif2}) is valid. Let us rewrite (\ref{phipsiif2}) in the form%
\begin{equation}
\phi_{x}-f_{0}^{2}\psi_{y}=0 \label{phipsi1}%
\end{equation}%
\begin{equation}
\phi_{y}+f_{0}^{2}\psi_{x}=0. \label{phipsi2}%
\end{equation}
From (\ref{phipsi2}) we have $f_{0}^{2}=-\phi_{y}/\psi_{x}$. Substituting this
expression into (\ref{phipsi1}) we obtain%
\[
\phi_{x}+\psi_{y}\frac{\phi_{y}}{\psi_{x}}=0
\]
or%
\[
\phi_{x}\psi_{x}+\phi_{y}\psi_{y}=0.
\]

\end{proof}

In general, solution of (\ref{phipsiif2}) or equivalently of the system
(\ref{phipsi1}), (\ref{phipsi2}) seems to be a difficult task. Nevertheless
for a quite general class of functions $f_{0}$ we can obtain solutions of
(\ref{phipsi1}), (\ref{phipsi2}) explicitly. Let us make the following suppositions.

\begin{condition}
\label{CondRho}Let $f_{0}$ be a function of some variable $\rho:$ $f_{0}%
=f_{0}(\rho)$ such that $\frac{\Delta\rho}{\left\vert \operatorname*{grad}%
\rho\right\vert ^{2}}$ is a function of $\rho$. We denote it by $s(\rho
)=\frac{\Delta\rho}{\left\vert \operatorname*{grad}\rho\right\vert ^{2}}$.
\end{condition}

The simplest example of such $\rho$ is of course any harmonic function.
Another important example is $\rho=\sqrt{x^{2}+y^{2}}$.

Consider the system (\ref{phipsi1}), (\ref{phipsi2}) and look for $\phi$ being
a function of $\rho:$ $\phi=\phi(\rho)$ (as we show below such solution always
exists). Then
\begin{equation}
\psi_{x}=-\frac{\rho_{y}}{f_{0}^{2}}\phi^{\prime},\qquad\psi_{y}=\frac
{\rho_{x}}{f_{0}^{2}}\phi^{\prime}. \label{psifromphi}%
\end{equation}
For the solubility of this system we obtain the following condition%
\[
\frac{\partial}{\partial x}(\frac{\rho_{x}}{f_{0}^{2}}\phi^{\prime}%
)+\frac{\partial}{\partial y}(\frac{\rho_{y}}{f_{0}^{2}}\phi^{\prime})=0
\]
which can be written as an ordinary differential equation%
\[
\phi^{\prime\prime}+(s-2\frac{f_{0}^{\prime}}{f_{0}})\phi^{\prime}=0.
\]
From here we have%
\[
\phi^{\prime}(\rho)=e^{-S(\rho)}f_{0}^{2}(\rho)
\]
where $S(\rho)=\int s(\rho)d\rho$.

With the aid of (\ref{psifromphi}) we can reconstruct $\psi$. Nevertheless we
are interested nor in $\phi$ neither in $\psi$ but in $\phi_{z}$ and $\psi
_{z}$ instead. Having them we construct the function $v=i(\phi_{z}F+\psi
_{z}G)=i(\phi_{z}/f_{0}+i\psi_{z}f_{0})$ which gives us a solution of
(\ref{Vek2}). We have
\begin{equation}
\phi_{z}=\phi^{\prime}\rho_{z}=e^{-S}f_{0}^{2}\rho_{z}\label{sol11}%
\end{equation}
and%
\begin{equation}
\psi_{z}=-\phi^{\prime}(\frac{\rho_{y}+i\rho_{x}}{f_{0}^{2}})=-ie^{-S}\rho
_{z}.\label{sol12}%
\end{equation}
Then we obtain the following solution of (\ref{Vek2})%
\[
v_{1}=i(\phi_{z}F+\psi_{z}G)=2if_{0}e^{-S}\rho_{z}.
\]
In a much the same way we can construct another solution of (\ref{Vek2})
looking for $\psi=\psi(\rho)$. Then
\[
\phi_{x}=f_{0}^{2}\rho_{y}\psi^{\prime},\qquad\phi_{y}=-f_{0}^{2}\rho_{x}%
\psi^{\prime}.
\]
and $\psi^{\prime}=e^{-S}/f_{0}^{2}$. Calculating $\phi_{z}$ and $\psi_{z}$ we
obtain%
\begin{equation}
\phi_{z}=ie^{-S}\rho_{z}\label{sol21}%
\end{equation}
and
\begin{equation}
\psi_{z}=\frac{e^{-S}}{f_{0}^{2}}\rho_{z}.\label{sol22}%
\end{equation}
Thus we arrive at the following solution of (\ref{Vek2})%
\[
v_{2}=i(\frac{ie^{-S}\rho_{z}}{f_{0}}+i\frac{e^{-S}}{f_{0}^{2}}\rho_{z}%
f_{0})=-2e^{-S}\frac{\rho_{z}}{f_{0}}.
\]

Denote
\begin{equation}
F_{I}=\frac{v_{1}}{2}\qquad\text{and\qquad}G_{I}=\frac{v_{2}}{2}. \label{F1G1}%
\end{equation}
Then $\operatorname*{Im}(\overline{F}_{I}G_{I})=e^{-2S}\left\vert
\operatorname*{grad}\rho\right\vert ^{2}>0$ and hence we have a generating
pair for (\ref{Vek2}) in explicit form. Note that $(F_{I},G_{I})$ is not a
successor of $(F,G)$ but a successor multiplied by $i$: $(F_{I},G_{I}%
)=i(F_{1},G_{1})$.

It is interesting to see what are the new solutions $f_{1}$ and $f_{2}$ of the
Schr\"{o}dinger equation (\ref{Schrc}) corresponding to $F_{I}$ and $G_{I}$.
As the procedure of converting solutions of (\ref{Vek2}) into solutions of
(\ref{Schrc}) requires on two steps to find a potential function by its
gradient, it will be more illustrative to consider here an example.

\begin{example}
\label{Ex1}Let $u(x,y)=x^{2}+y^{2}$. Then a particular solution of
(\ref{Schrc}) can be chosen in the form%
\[
f_{0}(x,y)=e^{xy}.
\]
Obviously $\rho(x,y)=xy$ being a harmonic function satisfies Condition
\ref{CondRho}. Then $F_{I}=ie^{xy}(y-ix)$ and $G_{I}=-e^{-xy}(y-ix).$
Returning to the notations of Section 3 (see the beginning of Section 5) we
have $F_{I}=\mathbf{k}e^{x_{1}x_{2}}(x_{1}-\mathbf{k}x_{2})=e^{x_{1}x_{2}%
}(x_{2}+\mathbf{k}x_{1})$ and $G_{I}=-e^{-x_{1}x_{2}}(x_{1}-\mathbf{k}x_{2})$.
Then the corresponding pair of solutions of (\ref{D+Mh}) has the form%
\[
\mathbf{F}_{1}=F_{I}\mathbf{j}=e^{x_{1}x_{2}}(-x_{1}\mathbf{i}+x_{2}%
\mathbf{j})\qquad\text{and\qquad}\mathbf{F}_{2}=G_{I}\mathbf{j}=-e^{-x_{1}%
x_{2}}(x_{2}\mathbf{i}+x_{1}\mathbf{j}).
\]
Next step (see Remark \ref{RemChain}) consists in finding the corresponding
functions $\Psi_{1}$ and $\Psi_{2}$ from equation (\ref{PsiFf}). Thus we
should reconstruct $\Psi_{1}$ and $\Psi_{2}$ from the equalities%
\[
\operatorname*{grad}\Psi_{1}=-x_{1}\mathbf{i}+x_{2}\mathbf{j}\qquad
\text{and\qquad}\operatorname*{grad}\Psi_{2}=-e^{-2x_{1}x_{2}}(x_{2}%
\mathbf{i}+x_{1}\mathbf{j}).
\]
Using the standard formula for finding the potential function from its
gradient we obtain%
\[
\Psi_{1}=-\frac{1}{2}(x_{1}^{2}-x_{2}^{2}-C_{1})\qquad\text{and\qquad}\Psi
_{2}=\frac{1}{2}(e^{-2x_{1}x_{2}}+C_{2})
\]
where $C_{1}$ and $C_{2}$ are arbitrary constants.

Now we can construct the vectors $\mathbf{g}_{1}=\operatorname*{grad}\Psi
_{1}/\Psi_{1}$ and $\mathbf{g}_{2}=\operatorname*{grad}\Psi_{2}/\Psi_{2}$:
\[
\mathbf{g}_{1}=\frac{2}{x_{1}^{2}-x_{2}^{2}-C_{1}}(x_{1}\mathbf{i}%
-x_{2}\mathbf{j)},\qquad\mathbf{g}_{2}=2(\frac{C_{2}}{e^{-2x_{1}x_{2}}+C_{2}%
}-1)(x_{2}\mathbf{i}+x_{1}\mathbf{j}).
\]
Noting that $\mathbf{h=}Df_{0}/f_{0}=x_{2}\mathbf{i}+x_{1}\mathbf{j}$ we
obtain two solutions $\mathbf{q}_{1}=\mathbf{h}+\mathbf{g}_{1}$ and
$\mathbf{q}_{2}=\mathbf{h}+\mathbf{g}_{2}$ for (\ref{Riccati}):%
\[
\mathbf{q}_{1}=(x_{2}+\frac{2x_{1}}{x_{1}^{2}-x_{2}^{2}-C_{1}})\mathbf{i}%
+(x_{1}-\frac{2x_{2}}{x_{1}^{2}-x_{2}^{2}-C_{1}})\mathbf{j}%
\]
and
\[
\mathbf{q}_{2}=(-x_{2}+\frac{2C_{2}x_{2}}{e^{-2x_{1}x_{2}}+C_{2}}%
)\mathbf{i}+(-x_{1}+\frac{2C_{2}x_{1}}{e^{-2x_{1}x_{2}}+C_{2}})\mathbf{j.}%
\]
Now we find the functions $\Phi_{1}$ and $\Phi_{2}$ which are solutions of the
equations $\operatorname*{grad}\Phi_{1,2}=\mathbf{q}_{1,2}$:%
\[
\Phi_{1}=\ln\left\vert x_{1}^{2}-x_{2}^{2}-C_{1}\right\vert +x_{1}x_{2}+C_{3}%
\]
and
\[
\Phi_{2}=x_{1}x_{2}+\ln\left\vert e^{-2x_{1}x_{2}}+C_{2}\right\vert +C_{4}.
\]
Then the corresponding solutions $f_{1}$ and $f_{2}$ of (\ref{Schrc}) have the
form%
\[
f_{1}=e^{\Phi_{1}}=d_{1}(x_{1}^{2}-x_{2}^{2}-C_{1})e^{x_{1}x_{2}}=d_{1}%
(y^{2}-x^{2}-C_{1})e^{xy}%
\]
and
\begin{equation}
f_{2}=e^{\Phi_{2}}=d_{2}(e^{-x_{1}x_{2}}+C_{2}e^{x_{1}x_{2}})=d_{2}%
(e^{-xy}+C_{2}e^{xy})\label{f2}%
\end{equation}
where $d_{1}$ and $d_{2}$ are arbitrary constant.

Thus starting with a particular solution of (\ref{Schrc}) we constructed two
classes of solutions for the same Schr\"{o}dinger equation.
\end{example}

\begin{example}
\label{ExCauchyTh} In order to illustrate the Cauchy integral theorem for the
Schr\"{o}dinger equation let us use the following two particular solutions
from the preceding example. Let $f_{0}(x,y)=e^{xy}$ and as $f_{1}$ we choose
the function from (\ref{f2}) when $C_{2}=0$ and $d_{2}=1$, $f_{1}%
(x,y)=e^{-xy}$. Both $f_{0}$ and $f_{1}$ are solutions of (\ref{Schrc}) with
the same potential $u$ in a whole plane. Thus we can apply Theorem
\ref{ThCauchySchr} and consider $\Gamma$ being for example a unitary
circumference with centre at the origin. Then
\[
\operatorname*{Re}\int_{\Gamma}\partial_{z}(\frac{f_{1}}{f_{0}}%
)dz+i\operatorname*{Im}\int_{\Gamma}f_{0}^{2}\partial_{z}(\frac{f_{1}}{f_{0}%
})dz
\]%
\[
=2\operatorname*{Re}\int_{\Gamma}(-y+ix)e^{-2xy}dz+2i\operatorname*{Im}%
\int_{\Gamma}(-y+ix)dz
\]%
\[
=2\operatorname*{Re}\int_{0}^{2\pi}i(-\sin\tau+i\cos\tau)e^{-2\cos\tau\sin
\tau}d\tau+2i\operatorname*{Im}\int_{0}^{2\pi}i(-\sin\tau+i\cos\tau)d\tau
\]%
\[
=-2\int_{0}^{2\pi}\cos\tau\cdot e^{-2\cos\tau\sin\tau}d\tau-2i\int_{0}^{2\pi
}\sin\tau d\tau.
\]
It is easy to see that both integrals are really equal to zero.
\end{example}

Let us calculate the characteristic coefficients for the pair $(F_{I},G_{I})$
defined by (\ref{F1G1}). We have
\[
a_{(F_{I},G_{I})}=0,\qquad b_{(F_{I},G_{I})}=-\frac{\partial_{z}f_{0}}{f_{0}%
},
\]

\[
A_{(F_{I},G_{I})}=\frac{\rho_{zz}}{\rho_{z}}-\frac{\rho_{z\overline{z}}}%
{\rho_{\overline{z}}},\qquad B_{(F_{I},G_{I})}=-\frac{f_{0}^{\prime}\rho
_{z}^{2}}{f_{0}\rho_{\overline{z}}}.
\]
Consider the equation
\[
\overset{1}{\phi}_{\overline{z}}F_{I}+\overset{1}{\psi}_{\overline{z}}G_{I}=0
\]
where $\overset{1}{\phi}$ and $\overset{1}{\psi}$ are real valued functions.
It has the form%
\begin{equation}
if_{0}^{2}\overset{1}{\phi}_{\overline{z}}-\overset{1}{\psi}_{\overline{z}}=0
\label{Kind2F1}%
\end{equation}
or as a system%
\begin{equation}
\overset{1}{\psi}_{x}+f_{0}^{2}\overset{1}{\phi}_{y}=0 \label{psiphi1}%
\end{equation}%
\begin{equation}
\overset{1}{\psi}_{y}-f_{0}^{2}\overset{1}{\phi}_{x}=0. \label{psiphi2}%
\end{equation}
Comparing with the system (\ref{phipsi1}), (\ref{phipsi2}) we note that
solutions of that system can be transformed into solutions of the system
(\ref{psiphi1}), (\ref{psiphi2}) in the following way%
\[
\overset{1}{\phi}=\psi\qquad\text{and\qquad}\overset{1}{\psi}=-\phi.
\]

\bigskip Let us calculate $F_{I}^{\ast}$ and $G_{I}^{\ast}$ using Definition
\ref{DefAdjoint}:%
\[
F_{I}^{\ast}=-\frac{f_{0}e^{S}}{\rho_{z}}\qquad\text{and\qquad}G_{I}^{\ast
}=-\frac{ie^{S}}{f_{0}\rho_{z}}.
\]
Consider the equation
\[
\overset{2}{\phi}_{\overline{z}}F_{I}^{\ast}+\overset{2}{\psi}_{\overline{z}%
}G_{I}^{\ast}=0
\]
where $\overset{2}{\phi}$ and $\overset{2}{\psi}$ are real valued functions.
It has the form%
\begin{equation}
if_{0}^{2}\overset{2}{\phi}_{\overline{z}}-\overset{2}{\psi}_{\overline{z}%
}=0.\label{kind2F1Adj}%
\end{equation}
Observe that it coincides with (\ref{Kind2F1}). Thus we obtain that the
function $v=\phi F_{I}+\psi G_{I}$ is $(F_{I},G_{I})$-pseudoanalytic iff the
function $W=\phi F_{I}^{\ast}+\psi G_{I}^{\ast}$ is $(F_{I}^{\ast},G_{I}%
^{\ast})$-pseudoanalytic.

Let us calculate
\[
B_{(F_{I}^{\ast},G_{I}^{\ast})}=\frac{\partial_{\overline{z}}f_{0}}{f_{0}}.
\]
That is the characteristic coefficient $b$ of a successor of $(F_{I}^{\ast
},G_{I}^{\ast})$ is equal to $-\partial_{\overline{z}}f_{0}/f_{0}=b_{(F,G)}$.
Thus, $(F,G)$ is a successor of $(F_{I}^{\ast},G_{I}^{\ast})$. This important
observation opens the way to obtain an infinite set of solutions of the
original Schr\"{o}dinger equation (\ref{Schrc}) if $f_{0}$ satisfies condition
\ref{CondRho}. Namely, we start with a solution of (\ref{Vek1}), for instance
with $F$. Then its $(F_{I}^{\ast},G_{I}^{\ast})$-antiderivative gives us an
$(F_{I}^{\ast},G_{I}^{\ast})$-pseudoanalytic function, more precisely the
corresponding functions $\phi$ and $\psi$ such that the function $W=\phi
F_{I}^{\ast}+\psi G_{I}^{\ast}$ is $(F_{I}^{\ast},G_{I}^{\ast})$%
-pseudoanalytic. Then we take these real valued functions $\phi$, $\psi$ and
consider the function $v=\phi F_{I}+\psi G_{I}$ which is $(F_{I},G_{I}%
)$-pseudoanalytic. That is $v$ satisfies (\ref{Vek2}) and hence the vector
$\mathbf{F}=v\mathbf{j}$ with the aid of the chain of actions described in
Remark \ref{RemChain} can be transformed into a solution of (\ref{Schrc}).
Taking the $(F,G)$-antiderivative of $iv$ we obtain another solution of
(\ref{Vek1}) and can start this cycle again. We can represent schematically
this procedure for obtaining an infinite sequence of solutions of (\ref{Vek2})
and consequently of (\ref{Schrc}) as the following diagram.
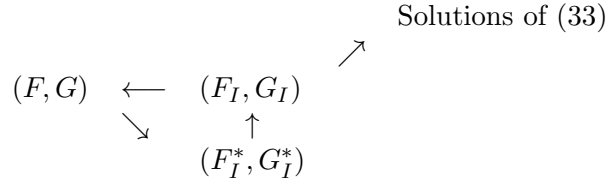
\[%
\begin{tabular}
[c]{lllll}
&  &  &  & Solutions of (\ref{Schrc})\\
&  &  & $\nearrow$ & \\
$(F,G)$ & $\longleftarrow$ & $(F_{I},G_{I})$ &  & \\
& $\searrow$ & \multicolumn{1}{c}{$\uparrow$} &  & \\
&  & $(F_{I}^{\ast},G_{I}^{\ast})$ &  &
\end{tabular}
\]
Let us consider how this procedure works on the following example.

\begin{example}
Here we use the same $u$ and $f_{0}$ as in Example \ref{Ex1}. Then $F=1/f_{0}%
$. We have that $F_{I}^{\ast}=-e^{xy}/(y-ix)$, $G_{I}^{\ast}=-ie^{-xy}%
/(y-ix)$. Consider%
\[
\ast\int_{0}^{z}Fd_{(F_{I}^{\ast},G_{I}^{\ast})}z=\operatorname{Re}\int
_{0}^{z}FG_{I}dz+i\operatorname{Re}\int_{0}^{z}FF_{I}dz
\]%
\[
=-\operatorname{Re}\int_{0}^{1}e^{-2xyt^{2}}%
(yt-ixt)(x+iy)dt+i\operatorname{Re}\int_{0}^{1}(xt+iyt)(x+iy)dt
\]%
\[
=\frac{e^{-2xy}-1}{2}+i\frac{(x^{2}-y^{2})}{2}.
\]
Thus we have $\phi=(e^{-2xy}-1)/2$ and $\psi=(x^{2}-y^{2})/2$. It is easy to
check that in fact $\phi_{\overline{z}}F_{I}^{\ast}+\psi_{\overline{z}}%
G_{I}^{\ast}=0$, so the function
\[
W=\phi F_{I}^{\ast}+\psi G_{I}^{\ast}=-\frac{1}{2(y-ix)}\left(  e^{-xy}%
-e^{xy}+ie^{-xy}(x^{2}-y^{2})\right)
\]
is $(F_{I}^{\ast},G_{I}^{\ast})$-pseudoanalytic. Now we can construct an
$(F_{I},G_{I})$-pseudoanalytic function as follows
\[
v=\phi F_{I}+\psi G_{I}=-\frac{(y-ix)}{2}\left(  e^{-xy}(x^{2}-y^{2}%
)-i(e^{-xy}-e^{xy})\right)  .
\]
Thus $v$ is a solution of (\ref{Vek2}) and applying the procedure described
above it can be transformed into a solution of (\ref{Schrc}).

Now multiplying the function $v$ by $i$ we obtain a solution of (\ref{WVek})
the $(F,G)$-antiderivative of which gives us a new solution of (\ref{Vek1})
and the cycle starts again.

Thus starting with a particular solution $F$ of (\ref{Vek1}) we obtain an
infinite sequence of solutions of (\ref{Schrc}).
\end{example}

\end{document}